# PYTHAGOREAN HARMONIC SUMMABILITY OF FOURIER SERIES

NASSAR H. S. HAIDAR

Euclidean Press LLC, 77 Ibrahim Abdul Al Str, Hamra, Beirut, Lebanon; email: nhaidar@suffolk.edu

**Abstract**. The paper explores the possibility for summing Fourier series nonlinearly via the Pythagorean harmonic mean. It reports on new results for this summability with introduction of new concepts like the smoothing operator and semi-harmonic summation. The smoothing operator is demonstrated to be Kalman filtering for linear summability, logistic processing for Pythagorean harmonic summability and linearized logistic processing for semi-harmonic summability. An emerging direct inapplicability of harmonic summability to seismic-like signals is shown to be resolvable by means of a regularizational asymptotic approach.

# 1. Introduction

It is well known that most $2L$-periodic $f(x) \in \mathcal{C}[0, 2L]$ can be represented by the Fourier series

$$f(x) := S(x) = \frac{a_0}{2} + \sum_{k=1}^{\infty}(a_k \cos k \tfrac{\pi}{L} x + b_k \sin k \tfrac{\pi}{L} x), \tag{1}$$

with the partial sums, see e.g. [1],

$$S_n(x) = \frac{a_0}{2} + \sum_{k=1}^{n}(a_k \cos k \tfrac{\pi}{L} x + b_k \sin k \tfrac{\pi}{L} x) = \sum_{k=0}^{n} \varphi_k(x), \tag{2}$$

given by the Dirichlet convolution integral, [2], [1],

$$S_n(x) = \frac{1}{2L} \int_{-L}^{L} D_n(\tau) f(x+\tau) \, d\tau, \tag{3}$$

with the Dirichlet kernel

$$D_n(x) = \frac{\sin(n+\tfrac{1}{2})\tfrac{\pi}{L} x}{\sin \tfrac{1}{2} \tfrac{\pi}{L} x}. \tag{4}$$

Moreover, in many applications, it is required to compare the smoothness of functional sums like $\sum_{k=0}^{n} \varphi_k(x)$ with $\sum_{k=0}^{n} \phi_n(k) \varphi_k(x)$, where $\phi_n(k)$ is a smooth function, which vanishes or decays for $k > n$. Such sums are called smoothed sums, [3]. This is what happens to take place with the Cesàro-Fejér sums, [1], [4],

$$\sigma_n(x) = \frac{a_0}{2} + \sum_{k=1}^{n} \phi_n(k)(a_k \cos k \tfrac{\pi}{L} x + b_k \sin k \tfrac{\pi}{L} x) = \frac{1}{(n+1)} \sum_{k=0}^{n} S_k(x), \tag{5}$$

associated with $S_n(x)$, where

$$\phi_n(k) = \left(1 - \frac{k}{n+1}\right) := \left\{\frac{n}{n+1}, \frac{n-1}{n+1}, \frac{n-2}{n+1}, \ldots, \frac{1}{n+1}\right\}. \tag{6}$$

---

Existence of $\sigma(x) = \lim\limits_{n \to \infty} \sigma_n(x)$, generated by $S_n(x)$, is known as $C_1$ summability of $S(x)$. Similarly, existence of $\hat{\sigma}(x)$, generated by $\sigma_n(x)$, is known as $C_2$ summability of $S(x)$, and so on.

Intuitively, one may expect a smoother behavior of the modified sums $\sum_{k=0}^{n} \phi_n(k)\, \varphi_k(x)$ because of the averaging involved (e.g. (5) for $\sigma_n(x)$) or because of the smoother cut-off at the end of the interval $[1,n]$ of summation, implying that those sums should be less susceptible to possible sudden violent oscillations in the size of the Fourier coefficients. This is visible analytically from the convolution integral expression which is now, [3],

$$\sigma_n(x) = \frac{1}{2L} \int_{-L}^{L} F_n(x) f(x+\tau)\, d\tau, \tag{7}$$

with the Fejér kernel

$$F_n(x) = \frac{1}{(n+1)} \left( \frac{\sin \frac{n+1}{2} \frac{\pi}{L} x}{\sin \frac{1}{2} \frac{\pi}{L} x} \right)^2. \tag{8}$$

The problem of establishing a criterion for convergence of Fourier series is a rather old one, and is still unresolvable in the form of a simple necessary and sufficient condition. For instance, bounded variation of $f(x)$ is sufficient but not necessary for convergence of $S(x)$. Continuity of $f(x)$ is neither necessary nor sufficient. Actually, there are $f(x)$ whose $S(x)$ converge at points of discontinuity and others whose $S(x)$ diverge at points of continuity. Classical summability of Fourier series was developed during the period 1897-1957, [1], as a complementary, or even "dual", property to their convergence. For instance, Fejér's summability theorem, [1-2] proves that at a point of continuity of $f(x)$, $S(x)$ is $C_1$ summable, so that continuity is at any rate sufficient for this summability.

It should be noted here, that all summation methods are based on averages. For example, in Cesàro-Fejér summation the average is arithmetical, while in the Abel-Poisson summation, the average is a harmonic function on the unit disk. Apart from these, the most important Fourier series summation methods are due to Riesz, Riemann, Bernstein-Rogosinski and de la Vallee-Poissin. Summation methods that are generated by a more-or-less arbitrary sequence of $\lambda$–multipliers, see e.g. [5], have also been studied.

Neoclassical summation methods for Fourier series were developed during the period 1960-1989 in [6-10], and in other works. Furthermore, during the period 1990-2017, contemporary summability research has ranged from factored Fourier series, [11], and product summability, [12], of these series, to generalizations for any orthogonal series, via sums based on Marcinkiewicz's $\Theta$–means, [13-14].

Despite the vastness of the surveyed theory for summing Fourier series, its practical applications continue to focus restrictively on:
i) Improving the representation of functions by Fourier series. For instance, if the $\sigma_n(x)$ sums, when $n \to \infty$, converge to $f(x)$ at its points of continuity, then they converge moreover uniformly on $[0, 2L]$ if $f(x) \in C[0, 2L]$. The partial sums $S_n(x)$ do not possess this property.
ii) Fourier-represented functional or perturbational, [15], analysis. For example, a function $f(x) := S(x)$ is essentially bounded iff $\exists$ a constant $M$ such that $|\sigma_n(x)| \leq M, \forall n$ & $x$.

On anoter note, a recent article, [16], of this author has reported on a minimal harmonic series for reconstructing an infinite Fourier series for $f(x) \in C[0, 2L]$. A series that is constructible by a minimal series interpolation (MSI) algorithm, [16]. The possibility for smoothing such series by linear (Cesàro-Fejér) summation has also been demonstrated in [16]. This paper is devoted to the subject of possible smoothing of Fourier series by nonlinear summations, particularly by Pythagorean harmonic and/or semi-harmonic summation.

The paper is organized as follows. After this introduction, section 2 introduces the concept of a symbolic smoothing operator and applies it to analysis of the Cesàro-Fejér $C_1$ linear summability.

The applications cover guarantees for smoothing by the $\mathcal{C}_1$ summability, that reveal certain affine and Kalman filtration features of the pertaining operator and establish conditions for its contraction mapping. Section 3 is devoted to nonlinear summabilities, with results on nonlinear processing by their smoothing operators and to their contraction mapping properties. Section 4 focuses on the Pythagorean harmonic sum and reports on a sharp result on its $\mathcal{I}_1$ summability. In section 5, an algorithm for a new semi-harmonic $\mathcal{J}_1$ summability is advanced. Unique linear processing and contraction mapping features are identified for the smoothing operator of this summability. Then it is demonstrated in section 6, how Pythagorean summabilities fail to apply to seismic-like (with zero mean) signals. Accordingly, a regularizational asymptotic approach is developed for handling the Pythagorean harmonic summability of such signals. The paper concludes by section 7.

## 2. Linear summation

It should be underlined, from the outset of this analysis, that the existing literature on summability of Fourier series is overwhelmingly immense, see e.g. [17-21]. Moreover, nonlinear summability, that is addressed in this work, should neither be confused with nonlinear Fourier analysis, as reviewed by Semmes in [17], nor with nonlinear convergence acceleration [22-23], [21].

Our starting point is accordingly an assertion that the $\sigma_n(x)$ sum of (5) can be expressed as

$$\sigma_n(x) = \frac{1}{(n+1)} \sum_{k=0}^{n-1} \frac{1}{2L} \int_{-L}^{L} D_k(\tau) f(x+\tau) d\tau + \frac{1}{(n+1)} S_n(x), \qquad (9)$$

which is the same as

$$\sigma_n(x) = W_\sigma(n,x) \;\boxed{\sigma}\; S_n(x), \qquad (10)$$

where the smoothing, of $S_n(x)$, operator

$$W_\sigma(n,x) = \frac{1}{(n+1)} \left\{ \sum_{k=0}^{n-1} \frac{1}{2L} \int_{-L}^{L} D_k(\tau) f(x+\tau) d\tau + (\cdot) \right\}, 0pt \qquad (11)$$

(which is the author's new concept) is symbolic (with a $\boxed{\sigma}$ multiplication) in the sense that it acts on $S_n(x)$ with the map:

$$W_\sigma(n,x) : S_n(x) \;\to\; \frac{1}{(n+1)} \left\{ \sum_{k=0}^{n-1} \frac{1}{2L} \int_{-L}^{L} D_k(\tau) f(x+\tau) d\tau + S_n(x) \right\}. 0pt \qquad (12)$$

The quantity to the right of $\to$ above is equal to $\sigma_n(x)$, as in (10).

Consideration of the substitutions

$$p_n = \frac{1}{(n+1)}, \qquad (13)$$

$$q_n(x) = \frac{1}{(n+1)} \sum_{k=0}^{n-1} \frac{1}{2L} \int_{-L}^{L} D_k(\tau) f(x+\tau) d\tau = \frac{1}{(n+1)} \sum_{k=0}^{n-1} S_k(x), \qquad (14)$$

in (12) illustrates, via

$$W_\sigma(n,x) \;\boxed{\sigma}\; S_n(x) = p_n S_n(x) + q_n(x) = \sigma_n(x), \qquad (15)$$

that $W_\sigma(n,x)$ is in fact an affine transformation of $S_n(x)$.

**Remark 1**. Making use of the identity
$q_n(x) = n\, p_n\, \sigma_{n-1}(x),$
in (15) transforms it to a recursive process:
$$W_\sigma(n,x) \;\boxed{\sigma}\; S_n(x) = \sigma_n(x) = n\, p_n\, \sigma_{n-1}(x) + p_n S_n(x), \qquad (16)$$

which is a statement that $\sigma_n(x)$ is conceivable as Kalman filtration (smoothing).

For any summation method to be applicable to divergent functional series, the method can satisfy, in a sufficient but not necessary sense, the following three properties, [3].

i) Regularity : The summation method should give correct answer for a convergent series.

ii) Linearity: a) If $\sum_{k=0}^{\infty} \varphi_k(x) = A(x)$ and $\sum_{k=0}^{\infty} \psi_k(x) = B(x)$, then $\sum_{k=0}^{\infty} [\varphi_k(x) + \psi_k(x)] = A(x) + B(x)$,

and (b) $\sum_{k=0}^{\infty} c\, \varphi_k(x) = c\, A(x), \forall c \in R$.

iii) Stability: If $\sum_{k=0}^{\infty} \varphi_k(x) = A(x)$, then $\sum_{k=0; k\neq j}^{\infty} \varphi_k(x) = A(x) - \varphi_j(x)$. (17)

In this respect, Cesàro-Fejér sum, $\sigma(x)$, is regular, linear and stable. But not every useful method for summing Fourier series can satisfy all the previous three requirements. Recently, a technique of nonlinear summation of power series was used in [24] for solving nonlinear evolution equations. It is our intention here to focus attention on similar alternative methods to explore their possible applicability and/or limitations for summing Fourier series.

For a study of the dynamics of $\sigma_n(x)$, with varying $n$, we shall suppress the x-variable by averaging over $[0, 2L]$, viz

$$\bar{\sigma}_n = \frac{1}{2L} \int_{-L}^{L} \sigma_n(x)\, dx, \; \bar{S}_n = \frac{1}{2L} \int_{-L}^{L} S_n(x)\, dx, \; \bar{\varphi}_n = \frac{1}{2L} \int_{-L}^{L} \varphi_n(x)\, dx; \bar{S}_0 = S_o, \bar{\varphi}_0 = \varphi_o = S_o, \quad (18)$$

to rewrite (16) as

$$\bar{\sigma}_n = \frac{n}{n+1} \bar{\sigma}_{n-1} + \frac{1}{n+1} \bar{S}_n = \bar{a}_n\, \bar{\sigma}_{n-1} + \bar{b}_n. \quad (19)$$

**Theorem 1**. *The difference equation (19) admits the number series solution*

$$\bar{\sigma}_n = \bar{\varphi}_0 + \sum_{k=1}^{n} \phi_n(k)\, \bar{\varphi}_k, \quad (20)$$

*where $\phi_n(k)$ is the smooth function (6).*

*Proof.* Equation (19) belongs in the class of first-order variable coefficient difference equations
$$y_n = a_n\, y_{n-1} + b_n, \; \forall n. \quad (21)$$
This class is analytically solvable in the form

$$y_n = \left(\prod_{i=1}^{n} a_i\right) y_0 + \sum_{k=1}^{n} \left(\prod_{i=k}^{n} \frac{a_i}{a_k}\right) b_k,$$

which is valid iff $a_n \neq 0, \; \forall n$.

In correspondence to (19), this solution can be shown to be

$$\bar{\sigma}_n = \left(\prod_{i=1}^{n} \frac{i}{i+1}\right) \bar{S}_0 + \sum_{k=1}^{n} \left(\prod_{i=k}^{n} \frac{i}{i+1}\right) \frac{\bar{S}_k}{k}, \quad (22)$$

and holds $\forall n \geq 0$.

Since $\prod_{i=1}^{n} \frac{i}{i+1} = \frac{1}{n+1}$, it is obvious that (22) is the same as

$$\bar{\sigma}_n = \frac{1}{(n+1)} \sum_{k=0}^{n} \bar{S}_k. \qquad (23)$$

Further substitution of $\bar{S}_k = \sum_{m=1}^{k} \bar{\varphi}_m$ in (23), with some algebra involving (6), lead to the required result. ∎

Furthermore, let $S_n(x) \in L^2[0, 2L]$, Hilbert space, endowed with the norm

$$\|S_n\| = \left\{ \int_{-L}^{L} S_n^2(x)\, dx \right\}^{1/2}, \qquad (24)$$

to state the results that follow.

**Theorem 2.** *The smoothing operator $W_\sigma(n,x)$ of a converging $S_n(x)$ is a contraction mapping when $n \in [N, \infty)$ with $N$ satisfying,*

$$N = \left( \frac{\|S_{N+1}\|}{\|S_N\|} - 1 \right)^{-1} - 1, \qquad (25)$$

*and is asymptotically cyclically stable.*

*Proof.* Consider

$$W_\sigma(n,x) \boxed{\sigma} S_{n+1}(x) = \frac{1}{n+2} S_{n+1}(x) + \frac{1}{n+2} \sum_{k=0}^{n} \frac{1}{2L} \int_{-L}^{L} D_k(\tau) f(x+\tau)\, d\tau = \sigma_{n+1}(x).$$

The summation above rewrites as

$$\frac{1}{n+2} \sum_{k=0}^{n-1} \frac{1}{2L} \int_{-L}^{L} D_k(\tau) f(x+\tau)\, d\tau + \frac{1}{n+2} \frac{1}{2L} \int_{-L}^{L} D_n(\tau) f(x+\tau)\, d\tau.$$

Therefore

$$\sigma_{n+1}(x) = \frac{n+1}{n+2} \frac{1}{n+1} \sum_{k=0}^{n-1} \frac{1}{2L} \int_{-L}^{L} D_k(\tau) f(x+\tau)\, d\tau + \frac{n+1}{n+2} \frac{1}{n+1} S_n(x) + \frac{1}{n+2} S_{n+1}(x)$$

$$= \frac{n+1}{n+2} \sigma_n(x) + \frac{1}{n+2} S_{n+1}(x). \qquad (26)$$

This relation can be represented as

$$W_\sigma(n,x) \boxed{\sigma} S_{n+1}(x) = G_\sigma(\sigma_n(x)) = \frac{n+1}{n+2} \sigma_n(x) + \frac{1}{n+2} S_{n+1}(x). \qquad (27)$$

Inductively, (27) regenerates (16) in the form

$$W_\sigma(n,x) \boxed{\sigma} S_n(x) = G_\sigma(\sigma_{n-1}(x)) = \frac{n}{n+1} \sigma_{n-1}(x) + \frac{1}{n+1} S_n(x). \qquad (28)$$

Take now the norm of the difference between (27) and (28):

$$\left\| W_\sigma \boxed{\sigma} S_{n+1} - W_\sigma \boxed{\sigma} S_n \right\| = \|G_\sigma(\sigma_n) - G_\sigma(\sigma_{n-1})\|$$

$$= \left\| \frac{n+1}{n+2} \left( \sigma_n - \frac{n(n+2)}{(n+1)^2} \sigma_{n-1} \right) + \frac{1}{n+2} S_{n+1} - \frac{1}{n+1} S_n \right\|. \qquad (29)$$

In view of the triangle inequality, (29) transforms to the inequality

$$\left\| W_\sigma \boxed{\sigma} S_{n+1} - W_\sigma \boxed{\sigma} S_n \right\| \leq \left\| \frac{n+1}{n+2} \left( \sigma_n - \frac{n(n+2)}{(n+1)^2} \sigma_{n-1} \right) \right\| + \left\| \frac{1}{n+2} S_{n+1} - \frac{1}{n+1} S_n \right\|. \qquad (30)$$

Convergence of $S_n(x)$ as $n \to \infty$ implies that for any $\delta \in (0,1]$, $\exists\, n = N(\delta) \in \mathbb{N}$, such that as of $n = N$, the following inequality

$$\left\| \frac{1}{n+2} S_{n+1} - \frac{1}{n+1} S_n \right\| \leq \delta, \tag{31}$$

is satisfied. Obviously $N$ is determined by solving

$$\left\| \frac{1}{N+2} S_{N+1} - \frac{1}{N+1} S_N \right\| = \delta. \tag{32}$$

Towards this end, we invoke the triangle inequality corollary

$$|\|\zeta\| - \|\xi\|| \leq \|\zeta - \xi\|, \tag{33}$$

with (32) to write

$$\left| \frac{1}{N+2} \|S_{N+1}\| - \frac{1}{N+1} \|S_N\| \right| \leq \delta = |\delta|.$$

This relation is the same as

$$\|S_{N+1}\| \leq \frac{N+2}{N+1} \|S_N\| + (N+2)\,\delta. \tag{34}$$

The earlier assumption on convergence of $S_n(x)$, allows for $\delta \to 0$ in (34), which yields the required formula, (25), for determination of $N$.

Obviously $\frac{\|S_{N+1}\|}{\|S_N\|} \geq 2$ leads, by (25), to $N = 0$, while $\frac{\|S_{N+1}\|}{\|S_N\|} = 1$ leads to $N = \infty$.

Furthermore, for $n \geq N$, $\frac{n(n+2)}{(n+1)^2} \to 1$ in (30) before $\frac{n+1}{n+2} \to 1$, which allows for

$$\left\| W_\sigma \boxed{\sigma} S_{n+1} - W_\sigma \boxed{\sigma} S_n \right\| = \|G_\sigma(\sigma_n) - G_\sigma(\sigma_{n-1})\| \leq \frac{n+1}{n+2} \|\sigma_n - \sigma_{n-1}\|, \tag{35}$$

i.e. $W_\sigma$ (or $G_\sigma$) is a contraction mapping. However, as $n \to \infty$, $\frac{n+1}{n+2} \to 1$, making $W_\sigma$ (or $G_\sigma$) asymptotically cyclic. Here the proof completes. ∎

Note incidentally that the contraction mapping property of $W_\sigma(n,x)$ should be a guarantee for the smoothing power of the associated $\sigma_n(x)$ sum.

## 3. Nonlinear summation

Actually, in addition to the arithmetic mean $\sigma_n(x)$ of (5), which is linear in elements of $\{S_k\}_{k=0}^n$, there are two other Pythagorean means, [25], that happen to be nonlinear in the previous elements. These are namely : the geometric mean

$$\gamma_n(x) = \left[ \prod_{k=0}^{n} S_k(x) \right]^{\frac{1}{n+1}}, \tag{36}$$

and the Pythagorean harmonic mean

$$\theta_n(x) = (n+1) \bigg/ \sum_{k=0}^{n} S_k^{-1}(x), \quad S_0(x) \neq 0 \; ; \tag{37}$$

not to be confused with the harmonic mean

$$\Omega(x) = \lim_{n \to \infty} \frac{1}{\log n} \sum_{k=0}^{n} \frac{S_{n-k}(x)}{(k+1)},$$

of the Riesz-Nörlund $\mathcal{H}_1$ summability, [26].

By the way, despite the above linearity of the $\sigma_n(x)$ summation, its pertaining smoothing operator $W_\sigma(n,x)$ is, by (15), nonlinear (affine). In actual fact, the same can also be said about most summations $\varsigma_n(x)$ which can be linked to a corresponding $\boxed{\varsigma}$ symbolic smoothing operator $W_\varsigma(n,x)$ viz

$$\varsigma_n(x) = W_\varsigma(n,x) \boxed{\varsigma} S_n(x), \quad \varsigma := \{\sigma, \gamma, \theta, \rho, \Omega, \ldots\}.$$

In particular, both Pythagorean mean symbolic operators

$$W_\gamma(n,x) = \left[\prod_{k=0}^{n-1} S_k(x)\right]^{\frac{1}{n+1}} (\cdot)^{\frac{1}{n+1}}, \tag{38}$$

and

$$W_\theta(n,x) = (n+1) \Big/ \left[\sum_{k=0}^{n-1} S_k^{-1}(x) + (\cdot)^{-1}\right], \tag{39}$$

act on the same $S_n(x)$, but with the distinct corresponding maps, viz

$$W_\gamma(n,x) : S_n(x) \to \left[\prod_{k=0}^{n-1} S_k(x)\right]^{\frac{1}{n+1}} S_n^{\frac{1}{n+1}}(x) = \gamma_n(x), \tag{40}$$

$$W_\theta(n,x) : S_n(x) \to (n+1) \Big/ \left[\sum_{k=0}^{n-1} S_k^{-1}(x) + S_n^{-1}(x)\right] = \theta_n(x). \tag{41}$$

**Proposition 1**. *The geometric summation $\gamma_n(x)$ is a power-type nonlinear smoothing process*

$$W_\gamma(n,x) \boxed{\gamma} S_n(x) = \gamma_n(x) = S_n^{\frac{1}{n+1}}(x)\, \gamma_{n-1}^{\frac{1}{n+1}}(x), \tag{42}$$

*which is asymptotically cyclic.*

*Proof.* Consideration of $\gamma_{n-1}(x) = \left[\prod_{k=0}^{n-1} S_k(x)\right]^{\frac{1}{n}}$ in (36) for $\gamma_n(x)$, together with the identity

$$\frac{1}{n+1} = \frac{1}{n} - \frac{1}{n(n+1)}, \tag{43}$$

yields (42). Obviously, $\gamma_n(x) = \gamma_{n-1}(x)$, as $n \to \infty$. ∎

**Proposition 2**. *The harmonic summation $\theta_n(x)$ is a fraction-type nonlinear smoothing process*

$$W_\theta(n,x) \boxed{\theta} S_n(x) = \theta_n(x) = \frac{\frac{1}{n} S_n(x) + S_n(x)}{\frac{1}{n} \theta_{n-1}(x) + S_n(x)} \theta_{n-1}(x), \tag{44}$$

*which is asymptotically cyclic.*

*Proof.* Consideration of $\theta_{n-1}(x) = n \Big/ \sum_{k=0}^{n-1} S_k^{-1}(x)$ in (41) for $\theta_n(x)$, together with some algebra ends up with (44). Also here, $\theta_n(x) = \theta_{n-1}(x)$, as $n \to \infty$. ∎

**Theorem 3**. *A sufficient condition for the smoothing operator $W_\gamma(n,x)$ to be a contraction mapping when n is large, but finite, is that*

$$\frac{\|\gamma_{n-1}\|}{\|S_n\|} > 1. \tag{45}$$

*Proof.* According to (42),

$$W_\gamma(n,x) \boxed{\gamma} S_{n+1}(x) = \gamma_n^{\frac{n+1}{n+2}}(x)\, S_{n+1}^{\frac{1}{n+2}}(x) = \gamma_{n+1}(x) = G_\gamma(\gamma_n(x)). \tag{46}$$

Take now the norm of the difference between (46) and (42):

$$\left\| W_\gamma \boxed{\gamma} S_{n+1} - W_\gamma \boxed{\gamma} S_n \right\| = \| G_\gamma(\gamma_n) - G_\gamma(\gamma_{n-1}) \|$$

$$= \left\| \gamma_n^{\frac{n+1}{n+2}} S_{n+1}^{\frac{1}{n+2}} - \gamma_{n-1}^{\frac{n}{n+1}} S_n^{\frac{1}{n+1}} \right\|$$

$$= \left\| \left( \frac{S_{n+1}}{\gamma_n} \right)^{\frac{1}{n+2}} \gamma_n - \left( \frac{S_n}{\gamma_{n-1}} \right)^{\frac{1}{n+1}} \gamma_{n-1} \right\|. \tag{47}$$

A necessary condition for the possibility of factoring $\| \gamma_n - \gamma_{n-1} \|$ out, when $n$ is large, but finite, of the relation for contraction mapping, is that

$$\frac{S_{n+1}(x)}{\gamma_n(x)} \approx \left( \frac{S_n(x)}{\gamma_{n-1}(x)} \right)^{\frac{n+2}{n+1}}. \tag{48}$$

In weak form, (42) is the same as

$$\frac{\|S_{n+1}\|}{\|\gamma_n\|} \approx \left( \frac{\|S_n\|}{\|\gamma_{n-1}\|} \right)^{\frac{n+2}{n+1}}. \tag{49}$$

Furthermore, $W_\gamma$ (or $G_\gamma$) can be a contraction mapping if (48), i.e. (49), is satisfied simultaneously with

$$\left( \frac{\|S_n\|}{\|\gamma_{n-1}\|} \right)^{\frac{1}{n+1}} = \sqrt[n+1]{\frac{\|S_n\|}{\|\gamma_{n-1}\|}} < 1. \tag{50}$$

This is clearly guaranteed when

$$\frac{\|S_n\|}{\|\gamma_{n-1}\|} < 1. \tag{51}$$

Notably, (51) is the same as the required (45), and (49) is always satisfied for large enough $n$. Here the proof ends. ∎

**Remark 2**. Condition (45) of Theorem 3 cannot hold when $n \to \infty$ because that contradicts with Proposition 1, on the asymptotic cyclic stability of $W_\gamma$.

# 4. Pythagorean harmonic summability

It should be noted here that for any sequence $(S_n) = (S_k)_{k=0}^n$, with $S_n > 0, \forall n$, the following famous inequality, [25],

$$\theta_n \leq \gamma_n \leq \sigma_n, \tag{52}$$

always holds. Incidentally, the inequality (52), which favors the harmonic mean, can be invoked to prove a lemma that follows.

**Lemma 1**. *Let $S_n(x) > 0, \forall n$. If $S(x)$ is $C_1$ summable to $\sigma(x)$, then it must also be $\mathcal{I}_1$ summable to $\theta(x)$. The opposite statement may, however, not be true.*

Furthermore, it is useful to underline, at this point, three features which add promise to employing the $\theta_n(x)$ harmonic mean to summability of $S(x)$.
i) The harmonic mean is often used in other disciplines, particularly in statistics, [25], to compute average ratios. Ratios that can incidentally be seen in the Fejér kernel (8) of the convolution integrals (7).
ii) This mean is less sensitive, than $\gamma_n(x)$ or $\sigma_n(x)$ to possible extremely large spikes in some $S_k(x)'s$ (like in point-wise divergence) of the $S(x)$ series.
iii) The apparent restriction of $S_0(x) \neq 0$ on direct application of $\theta_n(x)$ does not exclude, however, its

application in regularized form, to be shown later in subsection 6.1. Moreover, $\theta_n(x)$ has a more obvious relation to the summability of $S(x)$ than the alternative $\Omega$.

For these reasons, we shall focus attention, in the rest of this paper, on the $\theta_n(x)$ mean and, in the future, on a pertaining minimal, in the sense of [16], $\mathcal{I}_1$ summation series of $S(x)$. For an aim for this work, it perhaps suffices to merely establish relevance, robustness and /or possible limitations of this alternative new technique of summability.

**Theorem 4**. *A sufficient condition for the smoothing operator $W_\theta(n,x)$ to be a contraction mapping when n is large, but finite, is that*

$$\frac{\|\theta_{n-1}\|}{\|S_n\|} > 1. \tag{53}$$

*Proof.* Based on (44), let us write
$$\left\| W_\theta[\theta] S_{n+1} - W_\theta[\theta] S_n \right\| = \|G_\theta(\theta_n) - G_\theta(\theta_{n-1})\|$$
$$= \left\| \frac{\frac{1}{n+1}S_{n+1} + S_{n+1}}{\frac{1}{n+1}\theta_n + S_{n+1}} \theta_n - \frac{\frac{1}{n}S_n + S_n}{\frac{1}{n}\theta_{n-1} + S_n} \theta_{n-1} \right\|. \tag{54}$$

A necessary condition for the possibility of factoring $\|\theta_n - \theta_{n-1}\|$ out, when $n$ is large, but finite, of the relation for contraction mapping, is that

$$\frac{\frac{1}{n+1}S_{n+1}(x) + S_{n+1}(x)}{\frac{1}{n+1}\theta_n(x) + S_{n+1}(x)} \approx \frac{\frac{1}{n}S_n(x) + S_n(x)}{\frac{1}{n}\theta_{n-1}(x) + S_n(x)}, \tag{55}$$

which is the same as

$$(n+1)\frac{\theta_n(x)}{S_{n+1}(x)} - (n+2)\frac{\theta_{n-1}(x)}{S_n(x)} \approx -1. \tag{56}$$

In weak form, (56) is the same as

$$\left\| (n+1)\frac{\theta_n}{S_{n+1}} - (n+2)\frac{\theta_{n-1}}{S_n} \right\| \approx 1. \tag{57}$$

Then by the corollary (33), relation (57) converts to the inequality

$$(n+1)\frac{\|\theta_n\|}{\|S_{n+1}\|} - (n+2)\frac{\|\theta_{n-1}\|}{\|S_n\|} \leq 1. \tag{58}$$

$W_\theta$ (or $G_\theta$) can be a contraction mapping if (55), i.e. (58), is satisfied simultaneously with

$$\frac{\frac{1}{n}S_n(x) + S_n(x)}{\frac{1}{n}\theta_{n-1}(x) + S_n(x)} < 1. \tag{59}$$

Satisfaction of (53) in both fractions of (58), i.e.

$$\frac{\|\theta_n\|}{\|S_{n+1}\|} \approx \frac{\|\theta_{n-1}\|}{\|S_n\|} \approx C > 1, \tag{60}$$

transforms (58) to
$-C < 1$,
which is always true. Relation (60) also always guarantees satisfaction of (59). Here the proof completes. ∎

**Remark 3**. Condition (53) of Theorem 4 cannot hold when $n \to \infty$ because that contradicts with Proposition 2, on the asymptotic cyclic stability of $W_\theta$.

The nonlinear dynamics of (44), with varying $n$, cannot be freed of the $x$-dependence by an averaging process similar to (18). Nonetheless, for a representative fixed $x = x_0$, one can define
$\widetilde{\theta}_n = \theta_n(x_0)$, $\widetilde{S}_n = S_n(x_0)$ and $\widetilde{\varphi}_n = \varphi_n(x_0)$, $\tag{61}$

to analyze the, associated with (44), recurrence relation

$$\tilde{\theta}_n = \frac{\frac{1}{n}\tilde{S}_n + \tilde{S}_n}{\frac{1}{n}\tilde{\theta}_{n-1} + \tilde{S}_n} \tilde{\theta}_{n-1}. \tag{62}$$

This easily rewrites as

$$\tilde{\theta}_n = \frac{(n+1)}{n\left[1 + \frac{1}{n}\frac{\tilde{\theta}_{n-1}}{\tilde{S}_n}\right]} \tilde{\theta}_{n-1}. \tag{63}$$

Now, regardless of the magnitude and sign of $\frac{\tilde{\theta}_{n-1}}{\tilde{S}_n} \in R$, $\exists$ an $N \in \mathbb{N}$, after which

$$\frac{\tilde{\theta}_{n-1}}{n\,\tilde{S}_n} < 1, \forall n \geq N. \tag{64}$$

Consequently, for $n \geq N$, equation (63) tends to the generalized logistic equation

$$\tilde{\theta}_n = \frac{n+1}{n}\tilde{\theta}_{n-1} - \frac{n+1}{n^2\,\tilde{S}_n}\tilde{\theta}^2_{n-1}, \tag{65}$$

which is transformable, via the map

$$\tilde{\theta}_n = (n+1)\,\tilde{S}_n\,u_n, \tag{66}$$

to a standard logistic equation

$$u_n = \mathfrak{r}\,u_{n-1}(1 - u_{n-1}), \tag{67}$$

with a unit growth rate $\mathfrak{r}$.

The trajectories of (67) are known, see e.g. [27], to be generated by two fixed points (the first is at $u = 0$). Therefore one of the trajectories should necessarily satisfy $\lim_{n \to \infty} u_n = 0$. The other trajectory, namely due to $\mathfrak{r} = 1$, is fortunately, however, neither bifurcative nor chaotic. The same can, undoubtedly, be said about the trajectories of the precursor equation (65).

Despite the rather complicated analysis of the solution to (67), it is nonetheless, not impossible to attempt to linearize it. In this respect, relation (65) can readily be put into the form

$$\tilde{\theta}_n = \left[\frac{1}{n} + 1 - \frac{n+1}{n}\frac{1}{n}\frac{\tilde{\theta}_{n-1}}{\tilde{S}_n}\right]\tilde{\theta}_{n-1}. \tag{68}$$

In the special case of (64), when $\frac{\tilde{\theta}_{n-1}}{\tilde{S}_n} = \kappa > 1$, $1 - \kappa = -\alpha$, $\alpha > 0$, we have

$$\left[\frac{1}{n} + 1 - \frac{n+1}{n}\frac{1}{n}\frac{\tilde{\theta}_{n-1}}{\tilde{S}_n}\right] \approx \left(\frac{n-\alpha}{n} - \frac{\kappa}{n^2}\right). \tag{69}$$

Further consideration of (69) in (68) converts it to

$$\tilde{\theta}_n = \left(\frac{n-\alpha}{n} - \frac{\kappa}{n^2}\right)\tilde{\theta}_{n-1} + \tilde{e}_n, \tag{70}$$

where $\tilde{e}_n$ is a certain error induced by the approximation (69). In particular, if $\kappa = 2$, i.e. $\alpha = 1$, then

$$\tilde{\theta}_n = \left(\frac{n-1}{n} - \frac{2}{n^2}\right)\tilde{\theta}_{n-1} + \tilde{e}_n, \tag{71}$$

which is another kind of Kalman filtration, quite similar to (16), with its familiar smoothing power.

Moreover, one can always try to verify the existence of some reasonably representative trajectories for (62). Towards this end, let us rewrite (62) in the form

$$\tilde{\theta}_n = \frac{(n+1)}{1 + \sum_{k=0}^{n-1}\frac{\tilde{S}_n}{\tilde{S}_k}} \tilde{S}_n. \tag{72}$$

Since the number of $\tilde{\varphi}$–terms in $\tilde{S}_0$ is 1, in $\tilde{S}_1$ is 2, and in $\tilde{S}_n$ is $(n+1)$, we may consider a pair of extreme cases for such trajectories.

(i) If $\frac{\tilde{\varphi}_n}{\tilde{\varphi}_k} \approx 1, \forall\, n, k$, then

$$\frac{\widetilde{S}_n}{\widetilde{S}_0} = n+1, \quad \frac{\widetilde{S}_n}{\widetilde{S}_1} = \frac{n+1}{2}, \quad \frac{\widetilde{S}_n}{\widetilde{S}_2} = \frac{n+1}{3}, \ldots, \frac{\widetilde{S}_n}{\widetilde{S}_k} = \frac{n+1}{k+1},$$

and

$$\sum_{k=0}^{n-1} \frac{\widetilde{S}_n}{\widetilde{S}_k} = (n+1)\left[1 + \tfrac{1}{2} + \tfrac{1}{3} + \ldots + \tfrac{1}{n}\right] = (n+1)\ln(2n+1). \tag{73}$$

Eventually, for large enough $n$, we have

$$\widetilde{\theta}_n = \frac{\widetilde{S}_n}{\ln(2n+1)}, \tag{74}$$

indicating that $\lim_{n\to\infty} \widetilde{\theta}_n = 0$, $\forall \widetilde{S}_n$, a trajectory apparently associated with the first fixed point of (65).

(ii) If $\dfrac{\widetilde{S}_n}{\widetilde{S}_k} \approx 1, \forall\, n, k$, then

$$\widetilde{\theta}_n = \widetilde{S}_n, \tag{75}$$

is a trajectory generated by the second fixed point of (65), which is reasonably possible, despite its deficiency in any $\widetilde{S}_n$ smoothing features.

Next we report on a new basic result on the harmonic $\mathcal{I}_1$ summability of Fourier series. To simplify notation, use shall be made of the abbreviation

$$\sum\nolimits_p = \sum_{k=0}^{p} S_k^{-1}. \tag{76}$$

**Theorem 5**. *Let $S = S(x)$ be any bounded functional series over $I \subset R$. A sufficient condition for the $\mathcal{I}_1$ harmonic summability of this series is the existence of some finite $M > 0$ such that*

$$\left[S_{n+1}(x) - \sqrt{S_{n+1}^2(x) + 4M\,|S_{n+1}(x)|}\,\right] \leq 2\,\theta_n(x) \leq \left[S_{n+1}(x) + \sqrt{S_{n+1}^2(x) + 4M\,|S_{n+1}(x)|}\,\right], \tag{77}$$

*holds, $\forall n$ and $x \in I$.*

*Proof.* Constructively, invoke the harmonic mean $\theta_n(x)$ of (37) associated with the bounded sequence $(S_n) := (S_n(x))$. Consider next

$$\theta_{n+1} = (n+2) \Big/ \sum\nolimits_{n+1},$$

to define the increment

$$|\theta_n - \theta_{n+1}| = \left|\frac{(n+1)\sum_{n+1} - (n+2)\sum_n}{\sum_n \sum_{n+1}}\right|$$

$$= \left|\frac{(n+1)S_{n+1}^{-1} - \sum_n}{\sum_n \left[\sum_n + S_{n+1}^{-1}\right]}\right| = \left|\frac{\left[(n+1)S_{n+1}^{-1}/\sum_n\right] - 1}{\sum_n \left[\left(S_{n+1}^{-1}/\sum_n\right) + 1\right]}\right|.$$

In view of (37), this relation rewrites as

$$|\theta_n - \theta_{n+1}| = \left|\frac{\theta_n S_{n+1}^{-1} - 1}{S_{n+1}^{-1} + (n+1)/\theta_n}\right| = \left|\theta_n \frac{[\theta_n - S_{n+1}]}{[\theta_n + (n+1)S_{n+1}]}\right|.$$

Assume then that $m = n + r$ and consider

$$|\theta_n - \theta_m| = |(\theta_n - \theta_{n+1}) + (\theta_{n+1} - \theta_{n+2}) + \cdots + (\theta_{m-1} - \theta_m)|$$

$$\leq |\theta_n - \theta_{n+1}| + |\theta_{n+1} - \theta_{n+2}| + \cdots + |\theta_{m-1} - \theta_m| = \sum_{k=n}^{m-1}|\theta_{k+1} - \theta_{k+1}|. \tag{78}$$

Inductively, we may write
$$\left|\theta_{n+1} - \theta_{n+2}\right| = \left|\theta_{n+1} \frac{\theta_{n+1} - S_{n+2}}{\theta_{n+1} + (n+2) S_{n+2}}\right|.$$
Then the inequality (78) becomes
$$\left|\theta_n - \theta_m\right| \leq \left|\theta_n \frac{\theta_n - S_{n+1}}{\theta_n + (n+1) S_{n+1}}\right| + \left|\theta_{n+1} \frac{\theta_{n+1} - S_{n+2}}{\theta_{n+1} + (n+2) S_{n+2}}\right| + \cdots + \left|\theta_{m-1} \frac{\theta_{m-1} - S_m}{\theta_{m-1} + m S_m}\right|$$
$$= \sum_{k=n}^{m-1} \left|\theta_k \frac{\theta_k - S_{k+1}}{\theta_k + (k+1) S_{k+1}}\right|. \tag{79}$$

Due to boundedness of $S_n(x)$, with
$$\sup_{n;\, I} |S_n(x)| \leq K, \tag{80}$$
we may expect, for $k \geq m > n$ and $n$ large enough, $(k+1) S_{k+1} \gg \theta_k$. Therefore in the denominators of the fractions in (79) we may relatively ignore the $\theta'_n s$ to reinforce this inequality to
$$\left|\theta_n - \theta_m\right| \leq \sum_{k=n}^{m-1} \left|\theta_k \frac{\theta_k - S_{k+1}}{(k+1) S_{k+1}}\right|.$$

Replacement of all the integer factors $(k+1)$ in the numerators above by $(n+1)$, further strengthens this inequality to
$$\left|\theta_n - \theta_m\right| \leq \frac{1}{(n+1)} \sum_{k=n}^{m-1} \left|\theta_k \frac{\theta_k - S_{k+1}}{S_{k+1}}\right| = \frac{1}{(n+1)} \sum_{k=n}^{m-1} \left|\theta_k \left(\frac{\theta_k}{S_{k+1}} - 1\right)\right|$$
Now if
$$\sup_{n;\, I} \left|\theta_n(x) \left[\frac{\theta_n(x)}{S_{n+1}(x)} - 1\right]\right| \leq M, \tag{81}$$
when $S(x)$ is $\mathcal{I}_1$ summable, then we have
$$\theta_n^2 - S_{n+1} \theta_n \leq M \left|S_{n+1}\right|. \tag{82}$$
Clearly,
$$-2M \left|S_{n+1}\right| \leq \theta_n^2 - S_{n+1} \theta_n - M \left|S_{n+1}\right| \leq 0, \tag{83}$$
with the left constraint obviously redundant, leads to a range of admissible solutions
$$\theta_n^- \leq \theta_n \leq \theta_n^+, \tag{84}$$
where
$$\theta_n^\pm = \tfrac{1}{2}\left[\pm\sqrt{S_{n+1}^2 + 4M \left|S_{n+1}\right|} + S_{n+1}\right]. \tag{85}$$
Finally, substitute (37) in (83) - (84) to obtain (77). Then satisfaction of (81) subject to (82)-(83) leads to
$$\left|\theta_n - \theta_m\right| \leq \frac{rM}{(n+1)}.\text{0pt} \tag{86}$$
Therefore given $\varepsilon > 0$, if $n$ is chosen so large that $\frac{rM}{(n+1)} = \varepsilon$, $\forall r$ and $m = n + r$, then $(\theta_n)$ is a Cauchy sequence. Consequently, by the Cauchy criterion, see e.g. [2, 4], we infer that $(\theta_n)$ is a convergent sequence as $n \to \infty$. Conclusively, convergence of $(\theta_n)$ is obviously a sufficient condition for $S(x)$ to be $\mathcal{I}_1$ summable. ∎

**Claim 1**. The respective bounds $K$ and $M$ for $|S_n(x)|$ and $\left|\theta_n(x)\left[\frac{\theta_n(x)}{S_{n+1}(x)} - 1\right]\right|$ in Theorem 5 can possibly be equal.

*Proof.* Verification of this claim can directly be done by substituting $K$ for $M$ and for all the $S_k$'s of (77) to obtain

$(n + 1) = \frac{1}{2}\left[\sqrt{K^2 + 4K^2} + K\right]\frac{(n+1)}{K}.$

This turns out to lead to a remarkable correct fact :

$$\frac{1+\sqrt{5}}{2} > 1, \tag{87}$$

about the golden ratio $\varphi = \frac{1+\sqrt{5}}{2}$.

The alternative distinct $K$ and $M$ situation is certainly not ruled out. Here the proof ends. ∎

**Corollary 1**. *Let* $S = \sum_{k=0}^{\infty} v_k$ *be a convergent or divergent, but bounded, number series. A sufficient condition for the* $\mathcal{I}_1$ *summability of this series is the existence of some finite* $M > 0$ *such that the same relation (77) holds*, $\forall n$.

*Proof.* Same as the proof of the previous theorem with $K$ and $M$ standing for the bounds

$$|v_n| \leq K, \forall n, \text{ and } \left|\theta_n\left[\frac{\theta_n}{v_{n+1}} - 1\right]\right| \leq M, \forall n, \tag{88}$$

instead of the supremums, of (80) and (81), pertaining to $S_n = S_n(x)$.

**Example 1**. Consider the Grandi's series

$$S = \sum_{k=0}^{\infty}(-1)^k = 1 - 1 + 1 - 1 + 1 - \cdots$$

which is divergent but bounded. Its two accumulation points, 0 and 1, can respectively be identified by telescoping and bracketing.

The sequence of its partial sums

$(S_n) := 1, 0, 1, 0, 1, \cdots$

can further be employed with linear and nonlinear means to establish that $S$ is $\mathcal{C}_1$ summable to $\frac{1}{2}$, and $\mathcal{G}_1$ summable to 0.

Remarkably here $K = M = 1$; and when these are duly substituted in (77) we are guided again to $\frac{1+\sqrt{5}}{2} > 1$. This explains the $\mathcal{I}_1$ summability of this $S$ to 0, in agreement with its $\mathcal{G}_1$ summability.

Let us look back at $S(x)$ of (1) with the partial sums $S_n(x)$ of (2). These are representable [1, 2], in terms of the Dirichlet kernel (4) and Dirichlet integrals (3). Unlike $F_n(x)$, $D_n(x)$ is not $\geq 0$, but the two kernels are inter-related viz

$$F_n(x) = \frac{1}{(n+1)}\sum_{k=0}^{n} D_k(x). \tag{89}$$

Moreover, the harmonic mean $\theta_n(x)$ should satisfy

$$\theta_n(x) = (n+1)\left/\sum_{k=0}^{n}\left\{\frac{1}{2L}\int_{-L}^{L} D_k(\tau)f(x+\tau)d\tau\right\}^{-1}\right. \; ; \frac{a_0}{2} \neq 0, \tag{90}$$

which is apparently computationally expensive, due to the presence of the inverses of the convolution integrals in (90). This calls for a new constructive summation, to be called semi-harmonic, that turns out to be a linear processor and to simplify this problem.

# 5. A New Semi-Harmonic Summation

Assume the Fourier series S(x) not to be seismic-like, [28], i.e. its constant term $a_0/2 \neq 0$. Then invoke the representation (37) for $\theta_n(x)$ to employ in it a local averaging method to constructively create a semi-harmonic sum, $\rho_n(x)$, by the algorithm that follow. To simplify notation, use shall be made of the abbreviation $z = \mathfrak{M}\{C_k\}_{k=1}^n$ for the median element of an ordered set $\{C_k\}_{k=1}^n$.

**Algorithm 1**. Revisit (37) successively, starting with

$n = 1: \quad \theta_1 = \dfrac{2}{S_0^{-1} + S_1^{-1}} = \dfrac{2 S_0 S_1}{(S_0 + S_1)} = \dfrac{2}{Q} S_0 S_1$ ; use the arithmetic mean $\sigma_1$ of $Q$, viz $Q = z_1 2$

$\sigma_1, (z_1 = 1)$, to obtain $Q = 2\sigma_1$, then $\theta_1 = \dfrac{S_0}{\sigma_1} S_1 := \rho_1$.

$n = 2: \quad \theta_2 = \dfrac{3}{S_0^{-1} + S_1^{-1} + S_2^{-1}} = \dfrac{3 S_0 S_1 S_2}{S_0 S_1 + S_1 S_2 + S_2 S_0} = \dfrac{3}{Q} S_0 S_1 S_3$; represent $Q$, via $Q = (S_0 S_1 + S_1 S_2 + S_2 S_0) \approx z_2(S_0 + S_1 + S_2) = z_1 3 \sigma_2$, with

$z_2 = \mathfrak{M}\{S_0, S_1, S_2\} = S_1 \to \sigma_1$, for enhanced smoothing. (91)

This leads to $Q = \sigma_1 3 \sigma_2$, then $\theta_2 \approx \dfrac{S_0 S_1}{\sigma_1 \sigma_2} S_2 := \rho_2$.

$n = 3: \quad \theta_3 = \dfrac{4}{S_0^{-1} + S_1^{-1} + S_2^{-1} + S_3^{-1}} = \dfrac{3 S_0 S_1 S_2 S_3}{S_1 S_2 S_3 + S_0 S_2 S_3 + S_0 S_1 S_3 + S_0 S_1 S_2} = \dfrac{4}{Q} S_0 S_1 S_3 S_4$; represent $Q$, via $Q = (S_1 S_2 S_3 + S_0 S_2 S_3 + S_0 S_1 S_3 + S_0 S_1 S_2) \approx z_3(S_0 + S_1 + S_2 + S_4) = z_3 4 \sigma_3$, with

$z_3 = \mathfrak{M}\{S_1 S_3, S_2 S_3, S_0 S_1, S_0 S_2\} = S_1 S_2 \to \sigma_1 \sigma_2$, for enhanced smoothing. (92)

This leads to $Q = \sigma_1 \sigma_2 4 \sigma_3$, then $\theta_3 \approx \dfrac{S_0 S_1 S_2}{\sigma_1 \sigma_2 \sigma_3} S_3 := \rho_3$.

Continuation of this development, for any $n \geq 1$, allows for a generalization of the expression for $\rho_3(x)$ to

$$\rho_n(x) = \left( \prod_{k=0}^{n-1} \dfrac{S_k(x)}{\sigma_{k+1}(x)} \right) S_n(x). \tag{93}$$

Since $\sigma_{k+1}(x) = \dfrac{1}{(k+2)} \displaystyle\sum_{m=0}^{k+1} S_m(x)$, then

$$\rho_n(x) = \left( \prod_{k=0}^{n-1} (k+2) S_k(x) \Big/ \sum_{m=0}^{k+1} S_m(x) \right) S_n(x). \tag{94}$$

Existence of $\rho(x) = \displaystyle\lim_{n \to \infty} \rho_n(x)$, generated by $S_n(x)$, shall be called $\mathcal{J}_1$ summability of $S(x)$. It should be noted here that despite the identity $\rho_n(x) = \theta_n(x)$, only when $n = 1$, the $\mathcal{J}_1$ and $\mathcal{I}_1$ summabilities are structurally entirely different; $\rho_0(x)$ is undefined, while both sums are nonlinear in the elements of $\{S_k\}_{k=0}^n$. The successive local averagings, like (91)-(92), in the previous algorithm, suggest that the $\rho_n(x)$ sum should be smoother than $\theta_n(x)$, as to be verified later in this section.

As before, let us identify a smoothing operator $W_\rho(n,x)$ for $\rho_n(x) = W_\rho(n,x) \boxed{\rho} S_n(x)$, where

$$W_\rho(n,x) = \left( \prod_{k=0}^{n-1} (k+2) S_k(x) \Big/ \sum_{m=0}^{k+1} S_m(x) \right)(\cdot), \tag{95}$$

with $\boxed{\rho}$ = "·" distinctively, from $\boxed{\theta}$, is a direct multiplication,

$$W_\rho(n,x) \boxed{\rho} S_n(x) = W_\rho(n,x) S_n(x) 0pt, \tag{96}$$
i.e.
$$W_\rho(n,x) : S_n(x) \rightarrow \left( \prod_{k=0}^{n-1}(k+2) S_k(x) \bigg/ \sum_{m=0}^{k+1} S_m(x) \right) S_n(x) = \rho_n(x). \tag{97}$$

**Proposition 3**. *Let* $\sup_{n,I} |S_n(x)| < K < \infty$. *The semi-harmonic summation* $\rho_n(x)$ *is a linear smoothing process*

$$W_\rho(n,x) \boxed{\rho} S_n(x) = W_\rho(n,x) S_n(x) = \rho_n(x) = \left( (n+1) S_n(x) \bigg/ \sum_{m=0}^{n} S_m(x) \right) \rho_{n-1}(x), \tag{98}$$

*which is asymptotically cyclic.*

*Proof.* Consideration of $\rho_{n-1}(x) = \left( \prod_{k=0}^{n-2}(k+2) S_k(x) \bigg/ \sum_{m=0}^{k+1} S_m(x) \right) S_{n-1}(x)$ in (97) for $\rho_n(x)$, together with some algebra ends up with (98). Then replacement of all $S_n(x)$ by $K$, shows that $\rho_n(x) = \rho_{n-1}(x)$, as $n \rightarrow \infty$. ■

As with (44), the dynamics of the linear difference equation (98) can only be entertained when
$$\tilde{\rho}_n = \rho_n(x_0) \text{ and } \tilde{\sigma}_n = \sigma_n(x_0), \tag{99}$$
in the form
$$\tilde{\rho}_n = \left( (n+1) \tilde{S}_n \bigg/ \sum_{m=0}^{n} \tilde{S}_m \right) \tilde{\rho}_{n-1} = \frac{\tilde{S}_n}{\tilde{\sigma}_n} \tilde{\rho}_{n-1}. \tag{100}$$

Further substitutions of
$\tilde{\sigma}_n = \bar{a}_n \tilde{\sigma}_{n-1} + \bar{b}_n,$
and (4) in (100) lead to
$$\tilde{\rho}_n = \frac{(n+1) \tilde{S}_n}{n \tilde{\sigma}_{n-1} \left[ 1 + \frac{1}{n} \frac{\tilde{S}_n}{\tilde{\sigma}_{n-1}} \right]} \tilde{\rho}_{n-1}. \tag{101}$$

Now, as with (63), regardless of the magnitude and sign of $\frac{\tilde{S}_n}{\tilde{\sigma}_{n-1}} \in R$, $\exists$ an $M \in \mathbb{N}$, after which
$$\frac{\tilde{S}_n}{n \tilde{\sigma}_{n-1}} < 1, \forall n \geq M. \tag{102}$$
Clearly then, relation (101) tends, when $n \geq M$, to the linear difference equation
$$\tilde{\rho}_n = \left[ \frac{n+1}{n} \frac{\tilde{S}_n}{\tilde{\sigma}_{n-1}} - \frac{n+1}{n^2} \frac{\tilde{S}_n^2}{\tilde{\sigma}_{n-1}^2} \right] \tilde{\rho}_{n-1}. \tag{103}$$
Furthermore, if
$$\frac{\tilde{S}_n}{\tilde{\sigma}_{n-1}} \approx 1, \tag{104}$$
then
$$\tilde{\rho}_n = \left( \frac{n^2-1}{n^2} \right) \tilde{\rho}_{n-1} + \tilde{\epsilon}_n, \tag{105}$$
where $\tilde{\epsilon}_n$ is the error induced by adopting the approximation (104).

Clearly then, equation (105) is another kind of Kalman filtration, similar to ( but quite different from) (70). An indication that $\rho_n(x)$ is effectively a certain linearized version of $\theta_n(x)$.

**Theorem 6.** If $\dfrac{\|S_n\|}{\|S_{n+1}\|} > 1$ for some finite $n$, then the smoothing operator $W_\rho(n,x)$ is a contraction mapping and $\dfrac{\|\rho_{n-1}\|}{\|\rho_n\|} > 1$. Asymptotically, however, $\lim\limits_{n\to\infty} \dfrac{\|\rho_{n-1}\|}{\|\rho_n\|} = 1$.

*Proof.* According to (98),

$$W_\rho(n,x)\boxed{\rho}\, S_{n+1}(x) = \rho_{n+1}(x) = \left((n+2) S_{n+1}(x) \Big/ \sum_{m=0}^{n+1} S_m(x)\right) \rho_n(x) = G_\rho(\rho_n(x)). \tag{106}$$

Take now the norm of the difference between (106) and (98):

$$\left\| W_\rho\boxed{\rho} S_{n+1} - W_\rho\boxed{\rho} S_n \right\| = \|G_\rho(\rho_n) - G_\rho(\rho_{n-1})\|$$

$$= \left\| \left((n+2) S_{n+1} \Big/ \sum_{m=0}^{n+1} S_m\right) \rho_n - \left((n+1) S_n \Big/ \sum_{m=0}^{n} S_m\right) \rho_{n-1} \right\|.$$

A necessary condition for the possibility of factoring $\|\rho_n - \rho_{n-1}\|$ out, when $n$ is large, but finite, of the relation for contraction mapping, is that

$$\left((n+2) S_{n+1} \Big/ \sum_{m=0}^{n+1} S_m\right) \approx \left((n+1) S_n \Big/ \sum_{m=0}^{n} S_m\right), \tag{!07}$$

which is the same as

$$(n+2) S_{n+1} \sum_{m=0}^{n} S_m - (n+1) S_n \sum_{m=0}^{n+1} S_m \approx 0,$$

or

$$\left[\frac{n+2}{n+1}\frac{1}{S_n} - \frac{1}{S_{n+1}}\right] \sum_{m=0}^{n} S_m \approx 1. \tag{108}$$

In view of (98), $\sum_{m=0}^{n} S_m = (n+1) S_n \dfrac{\rho_{n-1}}{\rho_n}$, (108) becomes

$$\left[(n+2) - (n+1)\frac{S_n}{S_{n+1}}\right] \frac{\rho_{n-1}}{\rho_n} \approx 1. \tag{109}$$

In weak form, (109) is the same as

$$\left\| (n+2)\frac{\rho_{n-1}}{\rho_n} - (n+1)\frac{S_n}{S_{n+1}}\frac{\rho_{n-1}}{\rho_n} \right\| \approx 1. \tag{110}$$

Then by the corollary (37), relation (110) converts to the inequality

$$(n+2)\frac{\|\rho_{n-1}\|}{\|\rho_n\|} - (n+1)\frac{\|S_n \rho_{n-1}\|}{\|S_{n+1}\rho_n\|} \leq 1. \tag{111}$$

$W_\rho$ (or $G_\rho$) can be a contraction mapping if (107), i.e. (111), is satisfied simultaneously with

$$\left((n+1) S_n \Big/ \sum_{m=0}^{n} S_m\right) = \frac{\rho_n}{\rho_{n-1}} < 1, \text{ i.e. with } \frac{\rho_{n-1}}{\rho_n} > 1, \forall n < \infty, \,\&\, x \in I. \tag{112}$$

In weak sense, (112) is the same as

$$\frac{\|\rho_{n-1}\|}{\|\rho_n\|} \approx C > 1,$$

and its simultaneous satisfaction by (111) leads to

$$\frac{\|S_n\|}{\|S_{n+1}\|} \geq \frac{(n+2)}{(n+1)C}\left[C - \left(\frac{1}{n+2}\right)\right] > 1,$$

which would hold always and for all finite $n$. ∎

Hence for non-seismic like Fourier series $S(x)$, we may substitute (3) and (89) in (94) to obtain

$$\rho_n(x) = \left( \prod_{k=0}^{n-1} \frac{(k+2) \int_{-L}^{L} D_k(\tau) f(x+\tau) d\tau}{k+1 \sum_{m=0}^{k+1} \int_{-L}^{L} F_m(\tau) f(x+\tau) d\tau} \right) \frac{1}{2L} \int_{-L}^{L} D_n(\tau) f(x+\tau) d\tau, n \geq 1, \quad (113)$$

The main computational advantage of $\rho_n(x)$, via (113), over $\theta_n(x)$, is its apparent freedom of computing inverses of integrals that may tend to zero as $n \to \infty$.

## 6. Harmonic summability of seismic-like signals

Seismic-like (having $a_0 = 0$) $S(x)$ signals lead to $\theta_n(x) = \rho_n(x) = 0, \forall n$, which makes them redundant for $\mathcal{I}_1$ and $\mathcal{J}_1$ summation purposes. In this respect, a well known regularization method in the theory of ill-posed (unstable) problems of the first kind, due to Lavrentiev, [29], [30], consists in replacing the first kind problem with a parametrized problem of the second kind. Similar approximating considerations seem to be applicable for handling this situation. These lead us towards an asymptotic regularizational approach, that follows, for the $\mathcal{I}_1$ or $\mathcal{J}_1$ summabilities of seismic-like signals

**6.1.** *Asymptotic regularization procedure.* For $S(x)$ with $a_0 = 0$, the partial sums $S_n(x)$ are lead by $S_0(x) = 0$, which makes $\theta_n(x) = 0, \forall n$. Let us add then an arbitrary constant $\beta$, to $S(x)$ and obtain
$$Z(x) = S(x) + \beta := f^*(x) = f(x) + \beta. \quad (114)$$
The corresponding partial $Z_n(x)'s$ sums are tied with the the harmonic sum

$$\theta_n^*(x) = (n+1) \bigg/ \sum_{k=0}^{n} Z_k^{-1}(x), \ Z_0(x) = \beta, \ Z_k(x) = S_k(x) + \beta, \quad (115)$$

where

$$\theta_n^*(x) = (n+1) \bigg/ \sum_{k=0}^{n} \left\{ \frac{1}{2L} \int_{-L}^{L} D_k(\tau) f^*(x+\tau) d\tau \right\}^{-1}. \quad (116)$$

Now since $\lim_{n \to \infty} \theta_n^*(x) = Z(x)$, then

$$\left[ \lim_{n \to \infty} \theta_n^*(x) - \beta \right] = S(x), \forall \beta. \quad (117)$$

$\beta$ here happens to play a role similar, though different, from the regularization parameter, [29], [31], in the famous Tikhonov method of regularization. Nonetheless, arbitrariness of $\beta$ in the functional (117) demonstrates its regularizational nature. It is also remarkable how this procedure is also applicable to the semi-harmonic sum $\rho_n(x)$ that approximates $\theta_n(x)$. Indeed, for $S(x)$ with $a_0 = 0$,

$$\sigma_n = \frac{1}{(n+1)} \sum_{k=0}^{n} S_k = \frac{1}{(n+1)} \sum_{k=1}^{n} S_k. \quad (118)$$

But when $f^*(x) = S(x) + \beta$, $Z_0 = \beta$, with

$$\sigma_n^* = \frac{1}{(n+1)} \sum_{k=0}^{n} Z_k, \tag{119}$$

and

$$\rho_n^* = \beta \prod_{k=1}^{n} \frac{Z_k}{\sigma_k^*}. \tag{120}$$

Taking into consideration that

$$Z_k = S_k + \beta \text{ and } \sigma_k^* = \sigma_k^* + \frac{\beta}{(k+1)}, \tag{121}$$

in (92) reduces it to

$$\rho_n^* = (n+1)\beta \prod_{k=1}^{n} \frac{(S_k + \beta)}{[(k+1)\sigma_k^* + \beta]}. \tag{122}$$

This clearly allows for

$$\left[ \lim_{n \to \infty} \rho_n^*(x) - \beta \right] \approx S(x), \forall \beta, \tag{123}$$

which is approximately the same as (117).

Now after putting $\theta_n(x)$ (with $\theta_n^*(x)$) or even $\rho_n(x)$ (with $\rho_n^*(x)$) in a form of a regular summation series, like (5), it becomes possible to practically contemplate $\hat{\theta}_n(x)$ (with $\hat{\theta}_n^*(x)$) or even $\hat{\rho}_n(x)$ (with $\hat{\rho}_n^*(x)$) for minimal, in the sense of [16], summation series of $S(x)$.

# 7. Conclusion

In this paper the relevance of nonlinear averages to summing Fourier series in clearly demonstrated. The paper has been a platform for initiating research into the potential for the new subject of Pythagorean harmonic (and semi-harmonic) summability of Fourier series. The power of the new concept, of a symbolic smoothing, of $S_n(x)$, operator, in analyzing both linear and nonlinear summabilities has been demonstrated for the first time. This operator is demonstrated to be Kalman filtering for linear summability, logistic processing for Pythagorean harmonic summability and linearized logistic processing for semi-harmonic summability. An asymptotic regularizational technique is demonstrated to be applicable for summing up seismic-like signals. The newly reported semi-harmonic summability, with some unique smoothing properties, appears to exhibit much promise for future practical applications.

# References

**1**. N. K. Bary, *A Treatise on Trigonometric Series*, Pergamon Press, London, 1964.

**2**. E. C. Titchmarsh, *The Theory of Functions*, Oxford University Press, London, 1952.

**3**. H. Iwaniec, and E. Kowalski, *Analytic Number Theory*, AMS, Providence, 2004.

**4**. S. M. Nikolski, *A Course of Mathematical Analysis*, Vol. 2, Mir Publishers, Moscow, 1977.

**5**. A. Castillo, J. Chavez, and H. Kim, A Note on divergent Fourier series and $\lambda$-permutations, *The Australian Journal of Mathematical Analysis and Applications* 14(1), (2017), 1-9.

**6**. G. M. Petersen, Summability of a class of Fourier series, *Proceedings of the American Mathematical Society* 11(6), (1960), 994-998.

**7**. K. Kanno, On the absolute summability of Fourier series (II), *Tohoku Mathematical Journal* 13(12), (1961), 201-215.

**8**. F. C. Hsiang, On |C,1| summability factors of a Fourier series at a given point, *Pacific Journal of Mathematics* 33(1), (1970), 139-147.

**9**. P. Chandra, and G. D. Dikshit, On the |B| and |E,q| summability of a Fourier series, its conjugate series and their derived series, *Indian Journal of Pure and Applied Mathematics* 12(11), (1981), 1350-1360.

**10**. B. E. Rhoades, Matrix summability of Fourier series based on inclusion theorems II, *Journal of Mathematical Analysis & Applications* 130 (2), (1988), 525-537.

**11**. H. Bor, A note on local property of factored Fourier series, *Nonlinear Analysis, Theory, Methods & Applications* 64(3), (2006), 513-517.

**12**. H. K. Nigam, and K. Sharma, On (E,1)(C,1) summability of Fourier series and its conjugate series, *International Journal of Pure & Applied Mathematics* 82(3), (2013), 365-375.

**13**. F. Weisz, Θ-summability of Fourier series, *Acta Mathematica Hungarica* 103, (2004), 139-176.

**14**. F. Weisz, Marcinkiewicz summability of Fourier series, Lebesgue points and strong summability, *Acta Mathematica Hungarica* 153, (2017), 356-381.

**15**. A. H. Nayfe, *Perturbation Methods*, John Wiley & Sons, New York, 1973.

**16**. N. H. S. Haidar, Interpolatory minimal series for reconstructiong an infinite Fourier series, *Bulletin of Mathematical Analysis and Applications* 12(3), (2020), 34-45.
,
**17**. S. W. Semmes, Nonlinear Fourier analysis, *Bulletin of the American Mathematical Society* 20, (1989), 1-18.

**17**. K. Rauf, J. O. Omolehin, and D. J. Evans, Further results on strong summability of Fourier series, *International Journal of Computer Mathematics* 25, (2007), 331-339.

**19**. P. Padhy, U. Misra, and M. Misra, *Summability Methods and Its Applications*, LAP Lamhert, New Delhi, 2012.

**20**. M. Mursaleen, *Applied Summability Methods*, Springer Briefs in Mathematics, Springer, Berlin, 2014.

**21**. M. Mursaleen, and F. Basar, *Sequence Spaces, Topics in Modern Summability Theory*, Taylor & Francis Group, New York, 2020.

**22**. D. A. Smith, and W. F. Ford, Numerical comparisons of nonlinear convergence accelerators,

*Mathematics of Computation* 38(158), (1982), 481-499.

**23**. O. Costin, G. Luo, and S. Tanveer, Divergent expansion, Borel summability and 3D Navier-Stokes equation, *Philosophical Transactions of the Royal Society A* 366(1876), (2008), 2775-2788.

**24**. A. I. Zemlyanukhin, and A. V. Bochkarev, Nonlinear summation of power series and exact solutions of evolution equations, *Russian Mathematics* 62(1), (2018), 29-35.

**25**. Y. Qin, *Integral and Discrete Inequalities and Their Applications*, Vol I, Birkhäuser, Basel, 2016.

**26**. O. P. Varshney, On the absolute harmonic summability of a series related to Fourier series, *Proceedings of the American Mathematical Society* 10(5), (1959), 784-789.

**27**. R. M. May, Simple mathematical models with complicated dynamics, *Nature* 261, (1976), 459-467.

**28**. M. Bath, *Mathematical Aspects of Seismology*, Elsevier, Amsterdam, 1968.

**29**. V. P. Tanana, *Methods for Solution of Nonlinear Operator Equations*, VSP, Utrecht,1997.

**30**. N. H. S. Haidar, A Green's function approach to invertibility, *Mathematica Japonica* 50 (1), (1999), 10-18.

**31**. N. H. S. Haidar, The collocational double series inverse in quasi-linear regularizer form, *J. Inverse & Ill-Posed Problems* 7 (2), (1999), 127-144.